\documentclass[a4paper]{article}

\usepackage[utf8]{inputenc}
\usepackage[T1]{fontenc}

\usepackage[top=3cm, bottom=3cm, left=2.5cm, right=2.5cm]{geometry}
\usepackage{amsbsy, amssymb, amsmath, fancyhdr, enumerate, ae, aecompl, xfrac, xcolor}
\usepackage{hyperref}
\usepackage{graphicx}
\usepackage[labelformat=simple]{subcaption}   

\usepackage{tocloft}

\usepackage{makecell} 
\usepackage{multirow}
\usepackage{algorithm}
\usepackage{algpseudocode}

\usepackage{appendix}
\usepackage{tabto}
\usepackage{lscape}
\usepackage{pdfpages}

\DeclareMathOperator*{\argmin}{arg\,min}

\renewcommand\vec{\mathbf}

\author{Marcel Makovník\thanks{Corresponding author: marcel.makovnik@fmph.uniba.sk} \and Pavel Chalmovianský   \\  \small{Faculty of Mathematics, Physics and Informatics, Comenius University in Bratislava}}
\title{Curvature estimation for meshes via algebraic quadric fitting}
\date{April 2023}
\begin{document}

\maketitle
\begin{abstract}
We introduce the novel method for estimation of mean and Gaussian curvature and several related quantities for polygonal meshes. The algebraic quadric fitting curvature (AQFC) is based on local approximation of the mesh vertices and associated normals by a quadratic surface. The quadric is computed as an implicit surface, so it minimizes algebraic distances and normal deviations from the approximated point-normal neighbourhood of the processed vertex. Its mean and Gaussian curvature estimate is then obtained as the respective curvature of its orthogonal projection onto the fitted quadratic surface. Experimental results for both sampled parametric surfaces and arbitrary meshes are provided. The proposed method AQFC approaches the true curvatures of the reference smooth surfaces with increasing density of sampling, regardless of its regularity. It is resilient to irregular sampling of the mesh, compared to the contemporary curvature estimators. In the case of arbitrary meshes, obtained from scanning, AQFC provides robust curvature estimation.
\end{abstract}  
\section{Introduction}
Characterization of the shape is well described for smooth curves and surfaces. The basic notion of curvature is subject to differential geometry, whose rigorous review may be found in e.\ g.\ \cite{do2016differential} or \cite{spivak1999comprehensive}. With the development of computer graphics, the approximation of surfaces by piecewise linear surfaces -- polygon meshes -- became predominant. Hence, the natural question arises about how the geometric properties known for such surfaces (especially smooth) may be described for their piecewise linear approximations. The curvature estimation is widely used in areas like surface segmentation \cite{lee2015mesh, benzian20203d}, mesh simplification \cite{yao2015quadratic}, remeshing \cite{su2019curvature} or feature line extraction \cite{hildebrandt2005smooth} and even in architecture \cite{pottmann2007geometry}.

Approaches for the estimation of curvature may be divided into several categories. One of the most frequently used is the construction of discrete curvature operators based solely on the geometry of the mesh. The second approach estimates curvature by taking its value from the underlying approximating surface. 

We employ the approach, where we fit an implicit quadric to a small part of the mesh. In contrary to the majority of methods, our contribution is, that we use also vertex normals in the fitting and the quadric is fitted in algebraic sense. This means, we do not minimize Euclidean distance, which usually leads to solving least squares. Instead, we focus on algebraic distance which approximates the Euclidean distance well in the vicinity of the fitted quadric, see \cite{taubin1994distance}. Moreover, we minimize also the deviation between the normals associated to vertices and the gradient of the fitted quadric.  

The paper is organized as follows. In Section \ref{SEC-related-works}, we mention several frequently used approaches for approximation of curvature for meshes. Then, we provide the description of our approach for computation of curvatures using quadric fitting in Section \ref{SEC-curvatures}. The last part of our contribution, Section \ref{SEC-results}, is dedicated to experimental results, comparing our  approach to frequently used methods on publicly available data. The quality of the approximating mesh and its approximation of the curvature values are examined. The summary of our work is found in Section \ref{SEC-conclusion}. 

\section{Related works}\label{SEC-related-works}
The first works on approximation of curvatures on meshes were motivated not only by needs of geometric modelling, but also by many algorithms used in computer vision. One of the early approaches used construction of parametric surface, which is a graph of a function given by a bivariate polynomial. For each vertex of the triangle mesh, this polynomial approximated the distances between the vertices of one-ring and the tangent plane passing through that vertex \cite{hamann1993curvature}. One-ring neighbourhood of the inspected vertex was used also by Taubin in \cite{taubin1995estimating} to estimate the curvature tensor. The proposed algorithm works for triangle meshes, since the necessary vertex normals are computed from face normals. The discretization of the curvature theory using mesh paralellity was explored in \cite{bobenko2010curvature}.

One of the most frequently used approximations was introduced in \cite{meyer2003discrete}. The popularity of this method for meshes is caused by its simple implementation, since the mean and Gaussian curvatures are computed using the areas and angles of triangles adjacent to the explored vertex. In \cite{yang2006robust} the extraction of principle curvatures was achieved using principal component analysis of local neighbourhoods.

Another class of the estimation methods is to choose the curvature of the underlying surface. Semenova  estimated curvatures using quadric, which were represented as functions in Monge form \cite{semenova2004curvature}. The local fitting of a biquadratic Bézier patch, whose control net was obtained by least squares fitting, was utilized in \cite{razdan2005curvature}. Another method uses algebraic fit of spheres for either point sets or point-normal sets, which immediately gives the value of mean curvature by using the radius of the sphere \cite{guennebaud2007algebraic}. Another way of describing the underlying surface was presented in \cite{oztireli2009feature}, where moving least squares are utilized. One of the main advantages of using the underlying surface is, that it is independent on the graph structure of a polygon mesh, i.\ e.\ the curvature may be described also for point clouds. Estimation of the curvature per triangle  rather than vertex was suggested in \cite{zhihong2011curvature}. Each triangle, whose each corner was given by a vertex and an associated normal, was interpolated by triangular Bézier patch and the first and the second fundamental form was obtained from its parametrization.

A robust statistical approach for estimation of the curvature tensor was introduced in \cite{kalogerakis2007robust}. The provided algorithm is capable of correcting the noise in the meshes. Its extension for processing of point clouds was introduced in \cite{kalogerakis2009extracting}. Principled construction of discrete differential operators may be used for computation of curvatures of arbitrary polygon meshes, see \cite{de2020discrete}. Several analyses and empirical comparisons of curvature estimation methods are available in \cite{gatzke2004improved, gatzke2006estimating, digne2014numerical, vavsa2016mesh} or \cite{kronenberger2018empirical}.

\section{Curvature for implicit surfaces}\label{SEC-curvatures}
Since in our research we are interested in algebraic fitting of the surfaces, we enlist basic notions concerning implicit surfaces which are used throughout our work. Subsequently, we include the curvature related notions in terms of implicit surfaces. 

Let $ g: \mathbb{R}^3 \rightarrow \mathbb{R} $ be a differentiable real valued function of three variables. Then, an \textit{implicit surface} $ \mathcal{Z}(g) $ is the set of all zeros of $ g $, i.\ e.\
\begin{equation} \label{EQ-set-of-all-zeros}
\mathcal{Z}(g) = \{ (x, y, z)^\top \in \mathbb{R}^3 \: | \: g(x, y, z) = 0  \}.
\end{equation}
The gradient of the implicit function $ \nabla g(x, y, z) $ represents the normal vector at the point $ (x, y, z)^\top \in \mathcal{Z}(g) $. If the function $ q: \mathbb{R}^3 \rightarrow \mathbb{R} $ is given by
\begin{equation}
\begin{split}
q(x, y, z) &= a_{11}x^2 + a_{22}y^2 + a_{33}z^2 + a_{12}xy + a_{13}xz + a_{23}yz \\
 &+ a_{14}x + a_{24}y + a_{34}z + a_{44},
\end{split}
\end{equation}
with $ a_{11}, \ldots, a_{44} \in \mathbb{R} $ and at least one of the numbers  $ a_{11}, \ldots, a_{23} $ being non-zero, we obtain the surface $ \mathcal{Z}(q) $ known as a \textit{quadratic surface} or simply \textit{quadric}. 

Consider a regular surface $ \mathcal{S} \subset \mathbb{R}^3 $. For each point $ P \in \mathcal{S} $, we may determine the values $ \kappa_1 $ and $ \kappa_2 $, called \textit{principal curvatures}, that are the maximum and minimum of normal curvatures at that point. For parametrized surfaces, the principal curvatures $ \kappa_1 $ and $ \kappa_2 $ are easily computed using the coefficients of the \textit{first} and \textit{second fundamental forms} \cite{spivak1999comprehensive}. Subsequently, we may define the two well-known measures of curvature -- \textit{Gaussian curvature} $ K $ and \textit{mean curvature} $ H $ given by
\begin{equation}
K = \kappa_1 \kappa_2, \quad \text{and} \quad H = \frac{\kappa_1 + \kappa_2}{2}.
\end{equation}

For our needs, it is suitable to use curvature formulas, obtained from the function $ g: \mathbb{R}^3 \rightarrow \mathbb{R} $ defining the implicit surface $ \mathcal{Z}(g) $. More precisely, these formulas rely on the gradient $ \nabla g $ and Hessian $ \vec{H}_g $ of the function $ g $.  These formulas (and their proofs) for both implicit planar curves and implicit surfaces may be found in \cite{goldman2005curvature}. The principal curvatures are derived from the Gaussian and mean curvatures, as follows.

Let $ g: \mathbb{R}^3 \rightarrow \mathbb{R} $  and let $ V $ be the regular point of the corresponding surface $ \mathcal{Z}(g) $. The Gaussian curvature $ K $ at the point $ V $ is then given by
\begin{equation}\label{EQ-gaussian-curv}
K = \frac{\nabla g(V)^\top \, \text{adj}(\,\vec{H}_g(V)\,) \, \nabla g(V)}{\Vert \nabla g(V) \Vert ^4},
\end{equation}
and the mean curvature $ H $ at $ V $ is given by
\begin{equation}\label{EQ-mean-curv}
H = \frac{\nabla g(V)^\top \,\vec{H}_g(V) \,\nabla g(V) - \Vert \nabla g(V) \Vert^2 \,\text{tr}(\,\vec{H}_g(V)\,)}{2\Vert \nabla g(V) \Vert^3},
\end{equation}
where $ \nabla g(V) $ and $ \vec{H}_g(V) $ denote the value of the gradient and the Hessian of the function $ g $ at the point $ V = (x_V, y_V, z_V) \in \mathcal{Z}(g) $, respectively. 
For completeness, knowing the values $ K $ and $ H $ we can now easily compute the principal curvatures as
\begin{equation}\label{EQ-principal-curvatures}
\kappa_{1,2} = H \pm \sqrt{H^2 - K}.
\end{equation}

As introduced in \cite{koenderink1990solid}, one may define the measure of the size and shape using the polar coordinates in the $ (\kappa_1, \kappa_2) $-plane.
Let $ \kappa_1 $ and $ \kappa_2 $ be the principal curvatures of a point $ V $ of the regular surface $ \mathcal{S} \subset \mathbb{R}^3 $. Then, the number 
\begin{equation}
R = \sqrt{\frac{\kappa_1^2 + \kappa_2^2}{2}}
\end{equation}
is called \textit{curvedness} and the number
\begin{equation}
S = -\frac{2}{\pi}\arctan\frac{\kappa_1 + \kappa_2}{\kappa_1 - \kappa_2}
\end{equation}
is called \textit{shape index} at the point $ V $.

The interpretation of curvedness and shape index might be more intuitive in comparison to the traditional measures -- the mean and Gaussian curvature. Clearly,  $ R \geq 0 $ and $ -1 \leq S \leq 1 $. On contrary, the bounds of $ H $ and $ K $ vary from surface to surface. The curvedness describes, what is the difference between the surface and a plane in the vicinity of the point $ V \in \mathcal{S} $. The values of $ R $ approaching zero signify more gently curved points, while $ R = 0 $ implies that the point is planar (i.\ e.\ $ \kappa_1 = \kappa_2 = 0 $). The values of $ S $ determine the shape of a local patch at the point $ V $. This allows us to classify points -- whether they are elliptic, parabolic, hyperbolic or umbilic. Note, that the shape index is indeterminate for planar points. 

To compute curvedness and shape index for implicit surfaces, it is more efficient to omit the computation of principal curvatures and use the Gaussian and mean curvature instead. Then, the curvedness $ R $ and shape index $ S $ at the point $ V $ are given by
\begin{equation}\label{EQ-curvedness}
R = \sqrt{2H^2 - K},
\end{equation}
and
\begin{equation}\label{EQ-shape-index}
S = -\frac{2}{\pi}\arctan 2\sqrt{H^2 - K}.
\end{equation}

\section{Quadric fitting curvature}
In this section we propose the novel method, referred to as Algebraic quadric fitting curvature (AQFC), to determine mean and Gaussian curvature at vertex of a polygon mesh. The method AQFC is based on finding an underlying surface to a certain area of the mesh, similarly as \cite{guennebaud2007algebraic} and  \cite{oztireli2009feature}. More precisely, the surface is a quadric fitted to a certain set of points and their associated normals.

\subsection{Input data and quadric computation}

The proposed method may process any type of polygonal mesh $ \mathcal{M} = \left\{ \mathcal{V}_\mathcal{N}, \mathcal{E}, \mathcal{F} \right\} $, where $ \mathcal{V}_\mathcal{N} $ is the set of vertex-normals, i.\ e.\ each vertex is associated with a normal vector, and $ \mathcal{E} $ and $ \mathcal{F} $ are respectively the set of edges and faces. The vertex normals may be prescribed or we may use any of the methods of the computation listed in \cite{jin2005comparison}, e.\ g.\ simple averaging of adjacent face normals or averaging weighted by vertex angles.  

For each processed vertex-normal $ (V, \vec{n}) \in \mathcal{V}_\mathcal{N} $, we assign a neighbourhood $ \mathcal{H} $ by finding the smallest possible vertex ring (including $ (V, \vec{n}) $), containing fixed number of vertex-normals. Subsequently, we determine a quadric $ \mathcal{Z}(q) $ which is the best approximation of $ \mathcal{H} $ as follows. 

Let $ \mathcal{H} = \left\{ (V_i, \vec{n}_i) \right\}, \: i = 1, ..., m $. We define an objective function $ f: \mathbb{R}^{10} \rightarrow \mathbb{R} $ as
\begin{equation}\label{EQ-objective-function}
f(a_{11}, ..., a_{44}) = \sum_{i = 1}^m \left( \omega_{V_i} q(V_i)^2 + \omega_{\vec{n}_i} \| \vec{n}_i - \nabla q(V_i) \|^2 \right),
\end{equation}
where $ \omega_{V_i}, \omega_{\vec{n}_i} > 0 $ are fixed weights assigned to the vertices $ V_i $ and normals $ \vec{n}_i $, respectively. The arguments of the function $ f $ are the coefficients of the unknown quadric $ q $ and the function accumulates the weighted squared algebraic distances between the vertices and the surface $ \mathcal{Z}(q) $ and deviations between gradients of the quadric and input normal vectors of $ \mathcal{H} $.

The coefficients of the fitted quadric are obtained by finding such $ a_{11}, ..., a_{44} $ that the function $ f $ attains a minimal value. Since $ f $ is a quadratic function, its global minimum is located at a critical point of $ f $. Assuming the Hessian of $ f $ is positive definite, the coefficients are obtained by solving the system of linear equations 
\begin{equation}\label{EQ-system-partial}
\dfrac{\partial f}{\partial a_{ij}}(a_{11}, ..., a_{44}) = 0, \quad 1 \leq i \leq j \leq 4.
\end{equation}

If we used unit weights for both vertices and normals, the quadric would fit all the point-normals of the neighbourhood $ \mathcal{H} $ with the same importance. This may result in unexpected visual result, since we naturally expect, that the quadric should fit the vertices in the first place. Also, the quadric is assigned to the processed vertex-normal $ (V, \vec{n}) $. Hence, we may modify the quadric in such way, that the point-normals closer to the processed vertex-normal are approximated better, i.\ e.\ we use the weights
\begin{equation}\label{EQ-dd-weights}
\omega_{V_i} = \dfrac{1}{e^{\Vert V - V_i \Vert^4}},  \quad \omega_{\vec{n}_i} = \frac{10^{-4}}{e^{\Vert \frac{\vec{n}}{\Vert \vec{n} \Vert} - \frac{\vec{n}_i}{\Vert \vec{n}_i \Vert} \Vert^2}}, \quad i = 1, ..., m.
\end{equation}
These weights reflect the Euclidean distance between $ V_i $ and $ V $ and deviation between the normalized associated normal $ \vec{n}_i $ and the normalized processed normal $ \vec{n} $.

\subsection{Estimation of mean and Gaussian curvature}

\begin{figure}
\centering
\def\svgscale{1}
\begingroup%
  \makeatletter%
  \providecommand\color[2][]{%
    \errmessage{(Inkscape) Color is used for the text in Inkscape, but the package 'color.sty' is not loaded}%
    \renewcommand\color[2][]{}%
  }%
  \providecommand\transparent[1]{%
    \errmessage{(Inkscape) Transparency is used (non-zero) for the text in Inkscape, but the package 'transparent.sty' is not loaded}%
    \renewcommand\transparent[1]{}%
  }%
  \providecommand\rotatebox[2]{#2}%
  \newcommand*\fsize{\dimexpr\f@size pt\relax}%
  \newcommand*\lineheight[1]{\fontsize{\fsize}{#1\fsize}\selectfont}%
  \ifx\svgwidth\undefined%
    \setlength{\unitlength}{153.5993435bp}%
    \ifx\svgscale\undefined%
      \relax%
    \else%
      \setlength{\unitlength}{\unitlength * \real{\svgscale}}%
    \fi%
  \else%
    \setlength{\unitlength}{\svgwidth}%
  \fi%
  \global\let\svgwidth\undefined%
  \global\let\svgscale\undefined%
  \makeatother%
  \begin{picture}(1,0.76741405)%
    \lineheight{1}%
    \setlength\tabcolsep{0pt}%
    \put(0,0){\includegraphics[width=\unitlength,page=1]{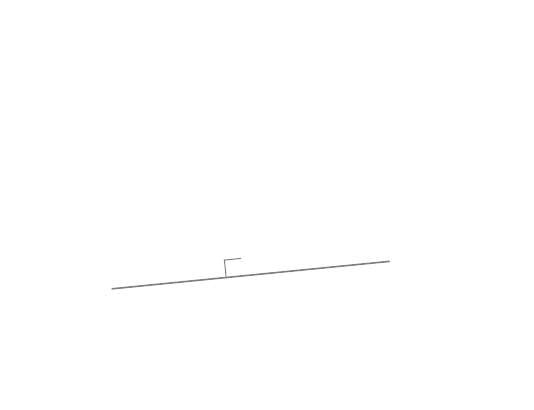}}%
    \put(0.42451734,0.65567837){\color[rgb]{0,0,0}\makebox(0,0)[lt]{\lineheight{1.25}\smash{\begin{tabular}[t]{l}$ V $\end{tabular}}}}%
    \put(0.0742965,0.01174154){\color[rgb]{0.49803922,0.49803922,0.49803922}\makebox(0,0)[lt]{\lineheight{1.25}\smash{\begin{tabular}[t]{l}$ \mathcal{Z}(q) $\end{tabular}}}}%
    \put(0,0){\includegraphics[width=\unitlength,page=2]{project-quadric.pdf}}%
    \put(0.46434698,0.17794227){\color[rgb]{0,0,0}\makebox(0,0)[lt]{\lineheight{1.25}\smash{\begin{tabular}[t]{l}$ p(V) $\end{tabular}}}}%
    \put(0,0){\includegraphics[width=\unitlength,page=3]{project-quadric.pdf}}%
    \put(0.45885751,0.43938358){\color[rgb]{0,0,0}\makebox(0,0)[lt]{\lineheight{1.25}\smash{\begin{tabular}[t]{l}$ \nabla q(p(V)) $\end{tabular}}}}%
  \end{picture}%
\endgroup%

\caption[Projection of a point onto an implicit surface.]{Projection of the point $ V $ onto the quadratic surface $ \mathcal{Z}(q) $. }\label{FIG-project-quadric}
\end{figure}

Let $ p(V) $ be the orthogonal projection of $ V $ onto the quadric $ \mathcal{Z}(q) $ given by 
\begin{equation}\label{EQ-projection-p}
p(V) = \argmin_{W \in \mathcal{Z}(q)} \Vert V - W \Vert. 
\end{equation}
Obviously, for the point of the quadratic surface $ V \in \mathcal{Z}(q) $, we have $ p(V) = V $.  Geometrically, the point $ p(V) $ is the footpoint of the point $ V $ onto the tangent plane of $ \mathcal{Z}(q) $ at the point $ p(V) $. In other words, the vector $ V - p(V) $ and the gradient $ \nabla q(p(V)) $ are linearly dependent, as seen in Figure \ref{FIG-project-quadric}.

Also, let $ s(\vec{n}, V) $ denote the gradient with corrected length given by
\begin{equation}
s(\vec{n}, V) = \| \vec{n} \| \frac{\nabla q(p(V))}{\| \nabla q(p(V)) \|}.
\end{equation}
The estimate of the mean curvature $ \tilde{H} $ and Gaussian curvature $ \tilde{K} $ assigned to the vertex-normal $ (V, \vec{n}) $ are computed as
\begin{equation}\label{EQ-gaussian-est}
\tilde{K} = \frac{s(\vec{n}, V)^\top \, \text{adj}(\,\vec{H}_q(p(V))\,) \, s(\vec{n}, V)}{\Vert s(\vec{n}, V) \Vert ^4}
\end{equation}
and
\begin{equation}\label{EQ-mean-est}
\hspace{-0.5em}
\tilde{H} = \frac{s(\vec{n}, V)^\top \,\vec{H}_q(p(V)) \,s(\vec{n}, V) - \Vert s(\vec{n}, V) \Vert^2 \,\text{tr}(\,\vec{H}_q(p(V))\,)}{2\Vert s(\vec{n}, V) \Vert^3}.
\end{equation}
Note, that we intentionally compute the curvatures at the point $ p(V) $ instead of the point $ V $. The reason is, that the point $ p(V) $ is the point of the quadric $ \mathcal{Z}(q) $, hence the respective gradient is not distorted by the displacement of the point $ V $. The possible error caused by the computation of the gradient at the point $ V $ would be reflected also in its norm. The displacement of $ V $ does not affect the Hessian, since $ q $ is given by quadratic polynomial. However, the error would be magnified in the computation of the mean and Gaussian curvatures, since the third and fourth power of the gradient norm is present in the denominators of the respective formulas. For the same reason, we rescale the original gradient $ \nabla q(p(V)) $ by the length of the vector $ \vec{n} $ in the estimations.

All remaining notions related to curvature -- principal curvatures, curvedness and shape index -- may be computed  using the respective formulas \eqref{EQ-principal-curvatures}, \eqref{EQ-curvedness} and \eqref{EQ-shape-index}.

\section{Experimental results}\label{SEC-results}
In this part, we explore the suitability of AQFC for estimating curvatures on polygon meshes.   Firstly, we describe the dataset, on which the experiments are made. For each mesh of the dataset we emphasize the features which may affect the precision of the curvature estimations. Then, we describe the measures used for the comparison. Finally, we provide the illustrations and measured values for the experiments and discuss the pros and cons of our method with respect to the selected competitors.

\subsection{Meshes for testing}\label{SUBSEC-meshes-for-testing}
We compare our method on several categories of polygonal meshes, whose description is provided below. The selected meshes are representative for the purposes of the error measurement, described in Subsection \ref{SUBSEC-error-measurement}.
\subsubsection{Regularly sampled surfaces}
The regular sampling of a surface is obtained from its parametric equations given by the function $ \vec{s}: \mathcal{D} \rightarrow \mathbb{R}^3 $, where $ \mathcal{D} \subset \mathbb{R}^2 $.

For testing, we consider the torus $ \vec{s}_{\mathbb{T}^2} $ with the minor radius equal to $ 1 $ and the major radius being $ 3 $. Its parametrization is given by
\begin{equation}\label{EQ-torus-parametrization}
\vec{s}_{\mathbb{T}^2}(\theta, \varphi) = 
\begin{pmatrix}
(3 + \cos \theta)\cos \varphi \\
(3 + \cos \theta)\sin \varphi \\
\sin \theta
\end{pmatrix}, \quad (\theta, \varphi) \in [0, 2\pi) \times [0, 2\pi).
\end{equation}
The mean and Gaussian curvature is computed directly from the parametrization as
\begin{equation}
H_{\mathbb{T}^2}(\theta, \varphi) = -\frac{3 + 2\cos \theta}{2(3 + \cos \theta)},
\end{equation}
and
\begin{equation}
K_{\mathbb{T}^2}(\theta, \varphi) = \frac{\cos \theta}{3 + \cos \theta},
\end{equation}
respectively. By computing extrema of functions $ H_{\mathbb{T}^2}(\theta, \varphi) $ and $ K_{\mathbb{T}^2}(\theta, \varphi) $, we get the bounds of the mean and Gaussian curvatures, which are 
\begin{equation}\label{EQ-torus-bounds}
\begin{split}
-0.65 &\leq H_{\mathbb{T}^2}(\theta, \varphi) \leq -0.25 \\
-0.5 &\leq K_{\mathbb{T}^2}(\theta, \varphi) \leq 0.25.
\end{split}
\end{equation}

By uniform sampling of the domain for parameters $ \theta $ and $ \varphi $, we constructed the quad meshes with 400, 3,600, 10,000,  360,000  and 1,000,000  vertices. Our aim is to inspect the behaviour of curvature estimation for increasing level of detail of the respective meshes.

\subsubsection{Surfaces with irregular sampling}  
Having a parametrization $ \vec{s}(\theta, \varphi) $ of the surface, we randomly choose the fixed number $ m $ of distinct pairs $ (\theta_i, \varphi_i) \in \mathcal{D} $, $ i = 1, ..., m $. The  samples are then obtained by plugging these random values of the parameters into the parametrization surface, i.\ e.\ the vertices $ V_i \in \mathcal{V} $ of the mesh are given by points $ \vec{s}(\theta_i, \varphi_i) $ for  $ i = 1, ..., m $. 

The faces of the mesh are given by the triangulation of the set of vertices $ \mathcal{V} $ obtained from Ball pivoting algorithm \cite{bernardini1999ball}. The input parameters for pivoting were manually adjusted to obtain the mesh with as few holes as possible. The remaining holes were enclosed, so the resulting mesh has no boundary.
For testing, we use the irregular sampling of the torus $ \vec{s}_{\mathbb{T}^2}(\theta, \varphi) $ \eqref{EQ-torus-parametrization} with 10,000 vertices. 

The mesh obtained by this procedure is irregular in terms of edge lengths, the ratio of obtuse and acute triangles, their areas and valence of the vertices. Our aim is to test the robustness of the curvature estimation with respect to such irregularities. 

Another instance of the irregular sampling is the one obtained from a unit sphere as follows. In the first step, we constructed the regular sampling of the sphere from its parametrization, similarly to the case of torus. Then, we split each quad face into two triangles by choosing a random diagonal of that quadrilateral.
 The multiple vertices at the poles of the sphere were merged to avoid degenerate faces. The resulting sphere has 
$ 482 $ vertices. 
The main feature of the unit sphere is, that its mean and Gaussian curvatures are equal and constant, that is 
\begin{equation}
H_{\mathbb{S}^2} = K_{\mathbb{S}^2} = 1.
\end{equation}

This sampling of the unit sphere is not as irregular as the sampling of the torus. Our goal is to check if even slight changes in the regularity of the sampling affect the precision of the estimation of the curvature.

\subsubsection{Meshes obtained from scanning}
We provide comparison also on the meshes -- Retheur and Damme Assise -- which are obtained by the surface reconstruction of the point cloud from the 3D scanning. These models are publicly available in the test dataset used in \cite{yifan2019patch}. The meshes have approximately 250,000 vertices.

\subsection{Tested methods and error measurement}\label{SUBSEC-error-measurement}
For comparison, we selected the curvature estimation by two commonly used state-of-the-art methods. The first one is DDGO \cite{meyer2003discrete}, which estimates curvature directly from the geometry of the input mesh. The second competitor is APSS \cite{guennebaud2007algebraic}, which is a method based on describing the underlying surface, similarly as our method AQFC. We used the Python library \textit{PyMeshLab} \cite{pymeshlab} for the computation of curvatures using DDGO and APSS.

For fairness, we consider all tested meshes without prescribed normals at vertices, even if they may be obtained from the parametrization or implicit equation of the sampled surface. For computation of vertex normals necessary for AQFC, we used the simple averaging of adjacent face normals.

For each method, we inspect separately on Gaussian curvature and mean curvature. Consider an input mesh $ \mathcal{M} = \{ \mathcal{V}, \mathcal{E}, \mathcal{F} \} $. For any vertex $ V \in \mathcal{V} $, let us denote by $ \tilde{K} $ the estimation of Gaussian curvature obtained by either AQFC, DDGO or APSS.  By $ K $ we denote the Gaussian curvature of the underlying surface. Then, we propose the following three measures $ K_{avg}, K_{min}, K_{max} \in \mathbb{R} $ given by
\begin{equation}\label{EQ-avg-gauss}
\begin{split}
K_{avg} &= \frac{1}{\vert \mathcal{V} \vert}\sum_{V \in \mathcal{V}} \vert \tilde{K} - K \vert, \\  
K_{min} &= \min_{V \in \mathcal{V}} \: \tilde{K}, \\
K_{max} &= \max_{V \in \mathcal{V}} \: \tilde{K}.
\end{split}
\end{equation}
The meshes in the dataset are obtained by sampling of a surface, where the curvature may be computed for each sample (vertex of the mesh). Hence, we may measure the absolute value of the difference between the estimated and true curvature at any mesh vertex. The value $ K_{avg} $ is a mean value of such differences and it provides us with the information how precisely each method approximates the curvatures of the underlying surface. Also, we inspect on the minimal and maximal values $ K_{min} $ and $ K_{max} $ to see, how well are the bounds of the curvature preserved. The difference between the bounds of estimated and true curvatures also shows the robustness of the estimates to outliers, caused by the features of the sampled mesh.

Analogously, we define the measures $ H_{avg}, H_{min}, H_{max} \in \mathbb{R} $ for mean curvature. These are dependent on the estimated values $ \tilde{H} $ and the true values $ H $ of mean curvature for any $ V \in \mathcal{V} $. Hence,
\begin{equation}\label{EQ-avg-mean}
\begin{split}
H_{avg} &= \frac{1}{\vert \mathcal{V} \vert}\sum_{V \in \mathcal{V}} \vert \tilde{H} - H \vert, \\  
H_{min} &= \min_{V \in \mathcal{V}} \: \tilde{H}, \\
H_{max} &= \max_{V \in \mathcal{V}} \: \tilde{H}.
\end{split}
\end{equation}

If the computed values of the proposed error measures for Gaussian and mean curvature are sufficiently small for certain method, we may assume that the method correctly estimates the principal curvatures, curvedness and shape index as well. Hence, the latter are not subject of our interest for the included experiments on regular and irregular sampling of the surfaces.

In the case of the meshes from 3D scans, we inspect on curvedness rather than Gaussian and mean curvature, since their bounds are unknown. We visualize the distribution of the curvedness across the tested meshes. We are interested in the presence of noise in the distribution. In this case, the noise is perceived as the difference in curvedness within small homogeneous  region of the mesh. For instance, in the areas with small curvedness (almost planar areas), we expect no vertices with high curvedness. 

\begin{figure}[ht]
\centering
\footnotesize
\begin{subfigure}{0.48\textwidth}
\centering
\def\svgwidth{\textwidth}
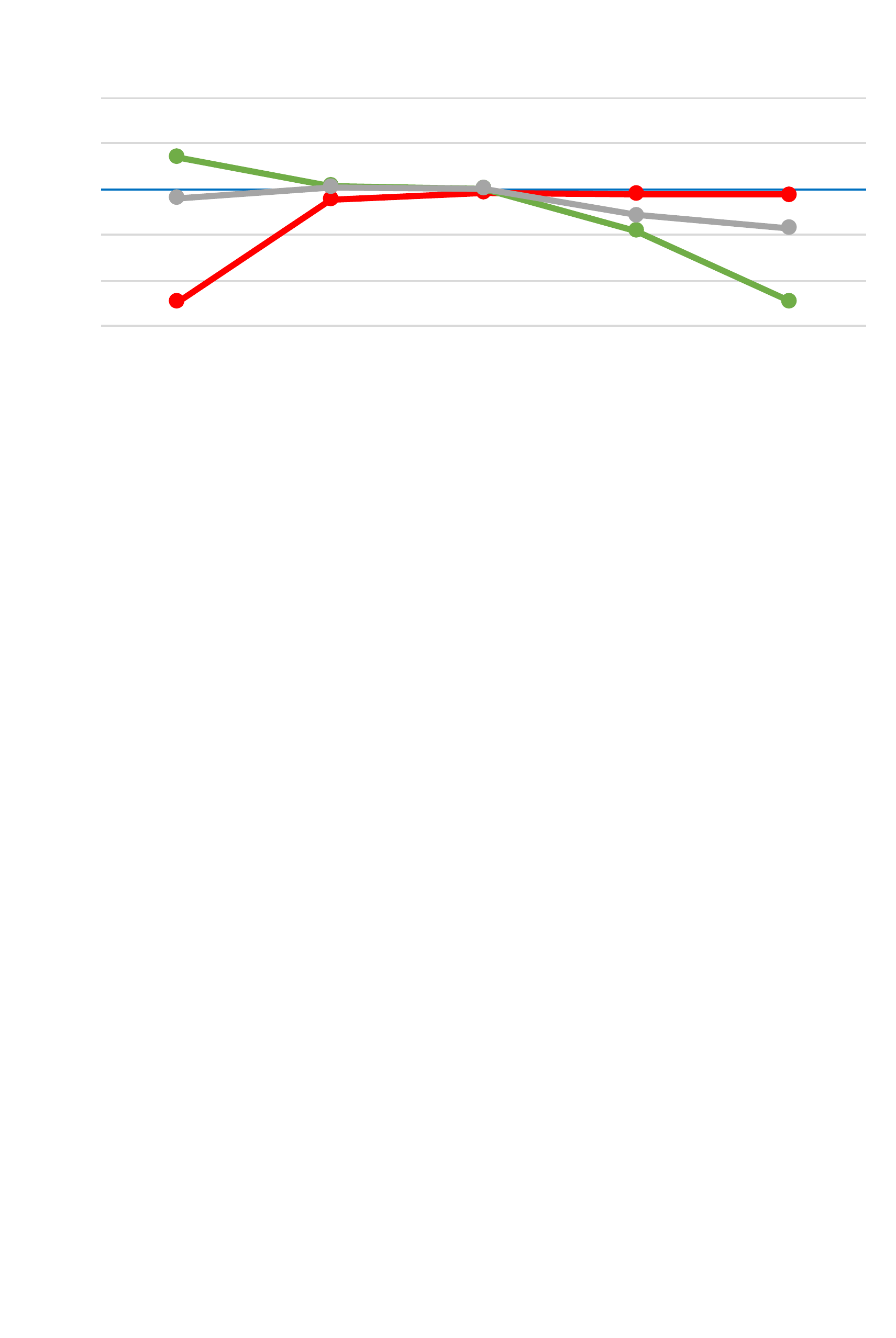
\caption{Mean curvature.}\label{FIG-mean_torus_reg}
\end{subfigure}\hfill
\begin{subfigure}{0.48\textwidth}
\centering
\def\svgwidth{\textwidth}
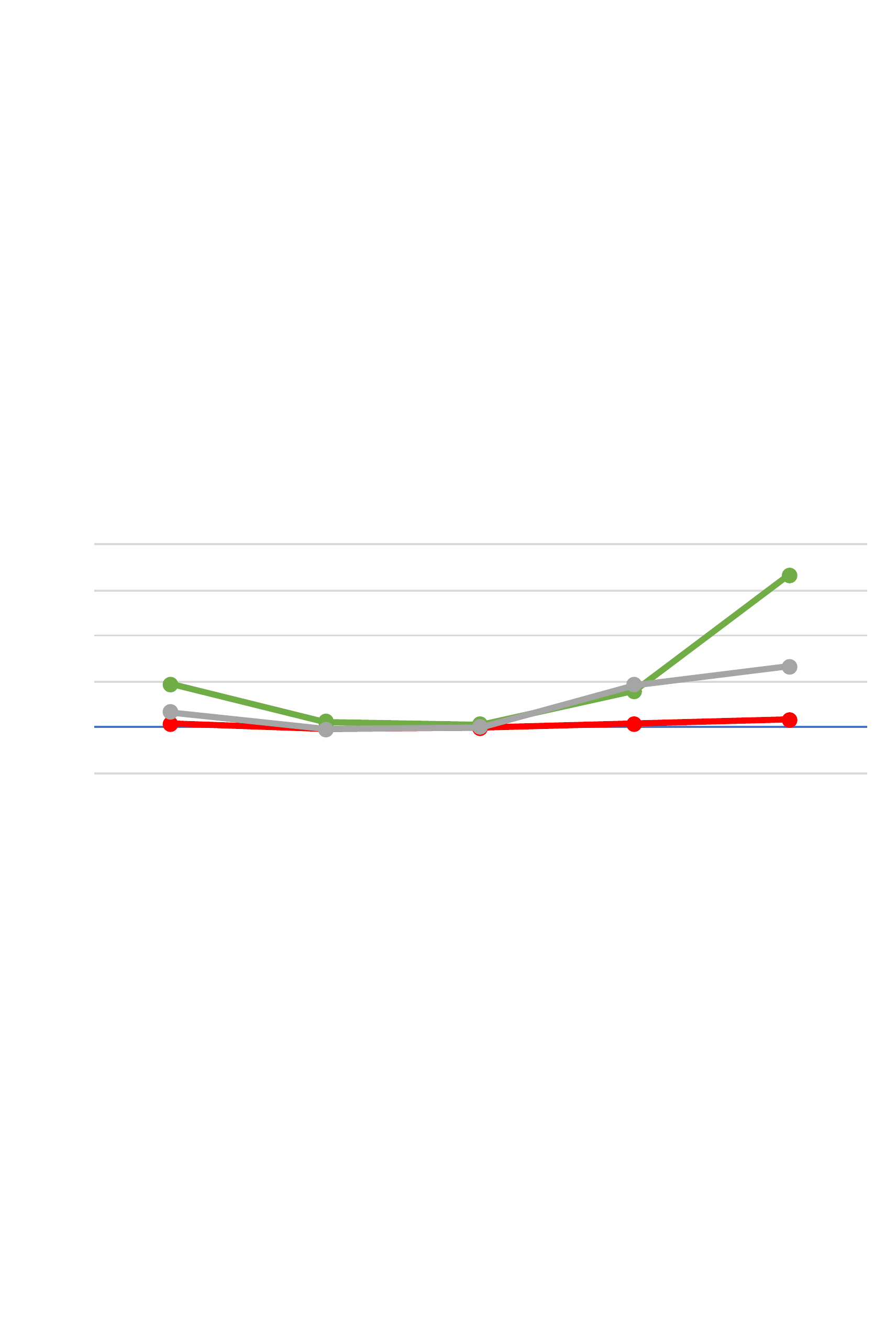
\caption{Gaussian curvature.}\label{FIG-gauss_torus_reg}
\end{subfigure}
\caption{The comparison of different methods of the minimal, maximal and average values for estimation of mean and Gaussian curvature on the regularly sampled torus at various levels of detail. The values on the horizontal axis denote the number of vertices of the sampled meshes.}\label{FIG-torus_reg}
\end{figure} 

\subsection{Comparison with other methods and discussion}
In this subsection, we provide the measurements of curvature on the meshes from Subsection \ref{SUBSEC-meshes-for-testing}, using the proposed method AQFC and methods DDGO and APSS. The figures of the meshes and their colourings were created using \textit{MeshLab} \cite{meshlab}.

\begin{figure}[ht]
\centering
\def\svgwidth{0.6\textwidth}
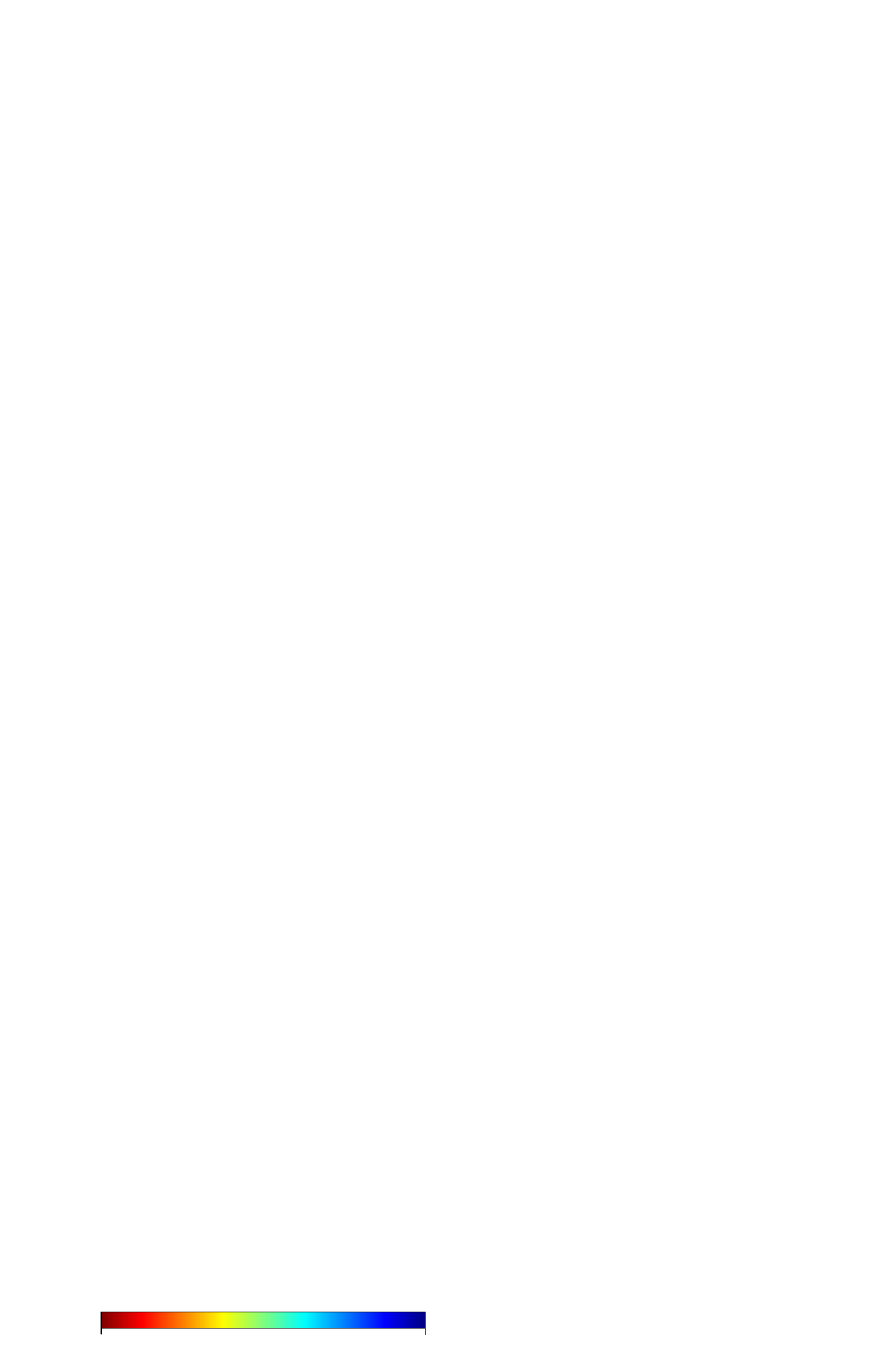
\caption[Estimation of curvature for irregularly sampled sphere using various methods.]{ The estimation of mean (left column) and Gaussian (right column) curvatures on the irregularly sampled sphere with $ 482 $ vertices.  } \label{FIG-curv-sphere-irreg}
\end{figure}

\subsubsection{Estimation precision with respect to level of detail}
For all methods, we measured minimal, maximal and average mean and Gaussian curvatures (see \eqref{EQ-avg-gauss} and \eqref{EQ-avg-mean}) of the regularly sampled torus at various level of density, as described in Subsection~\ref{SUBSEC-meshes-for-testing}. The values of the measurements are provided in Figure~\ref{FIG-mean_torus_reg} and \ref{FIG-gauss_torus_reg}, respectively.

\begin{figure}[ht]
\small
\centering
\begin{subfigure}{0.48\textwidth}
\def\svgwidth{\textwidth}
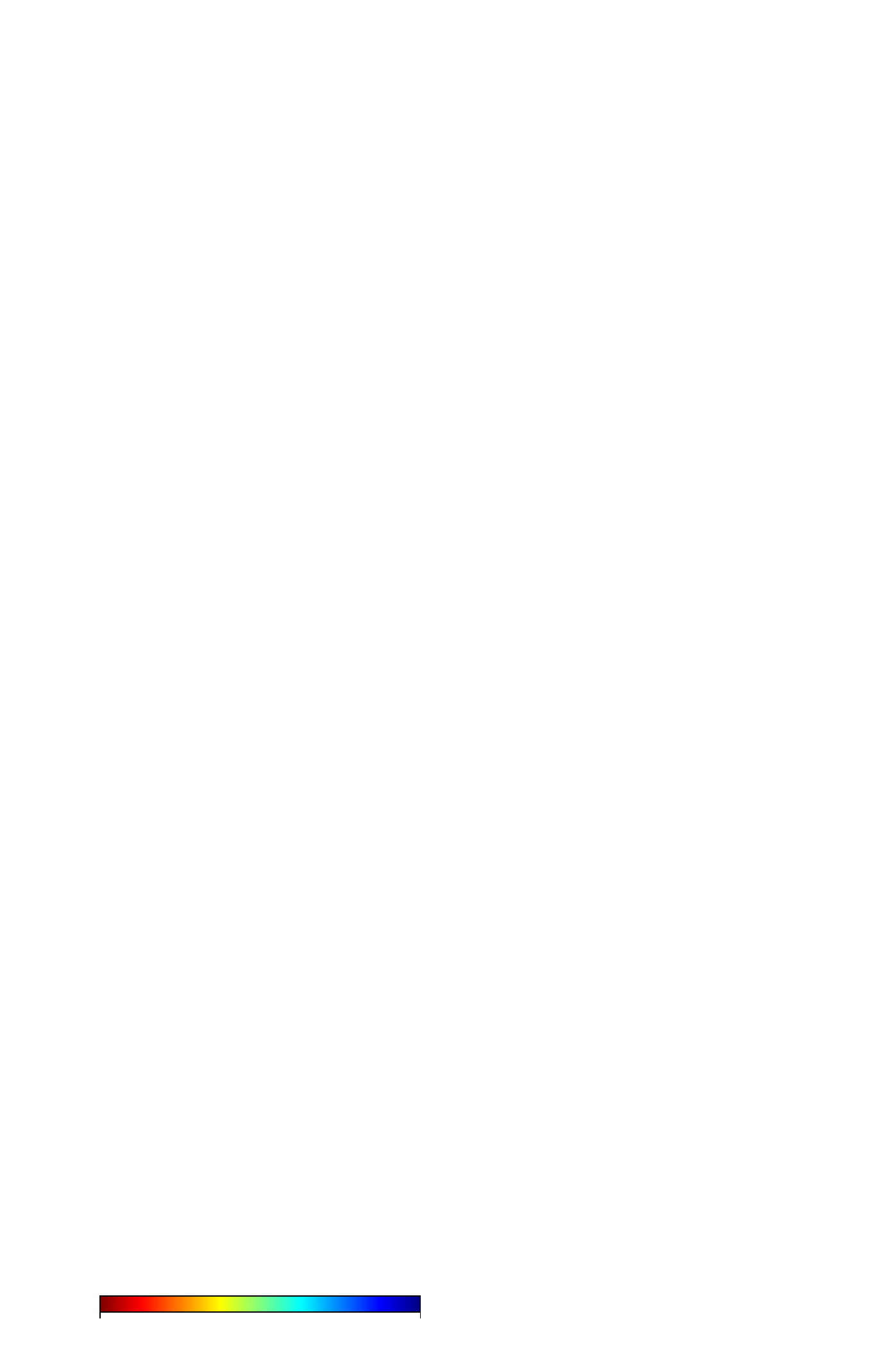
\caption[Estimation of curvature for regularly sampled torus using different methods.]{Regular sampling.} \label{FIG-curv-torus-reg}
\end{subfigure}\hfill
\begin{subfigure}{0.48\textwidth}
\centering
\def\svgwidth{\textwidth}
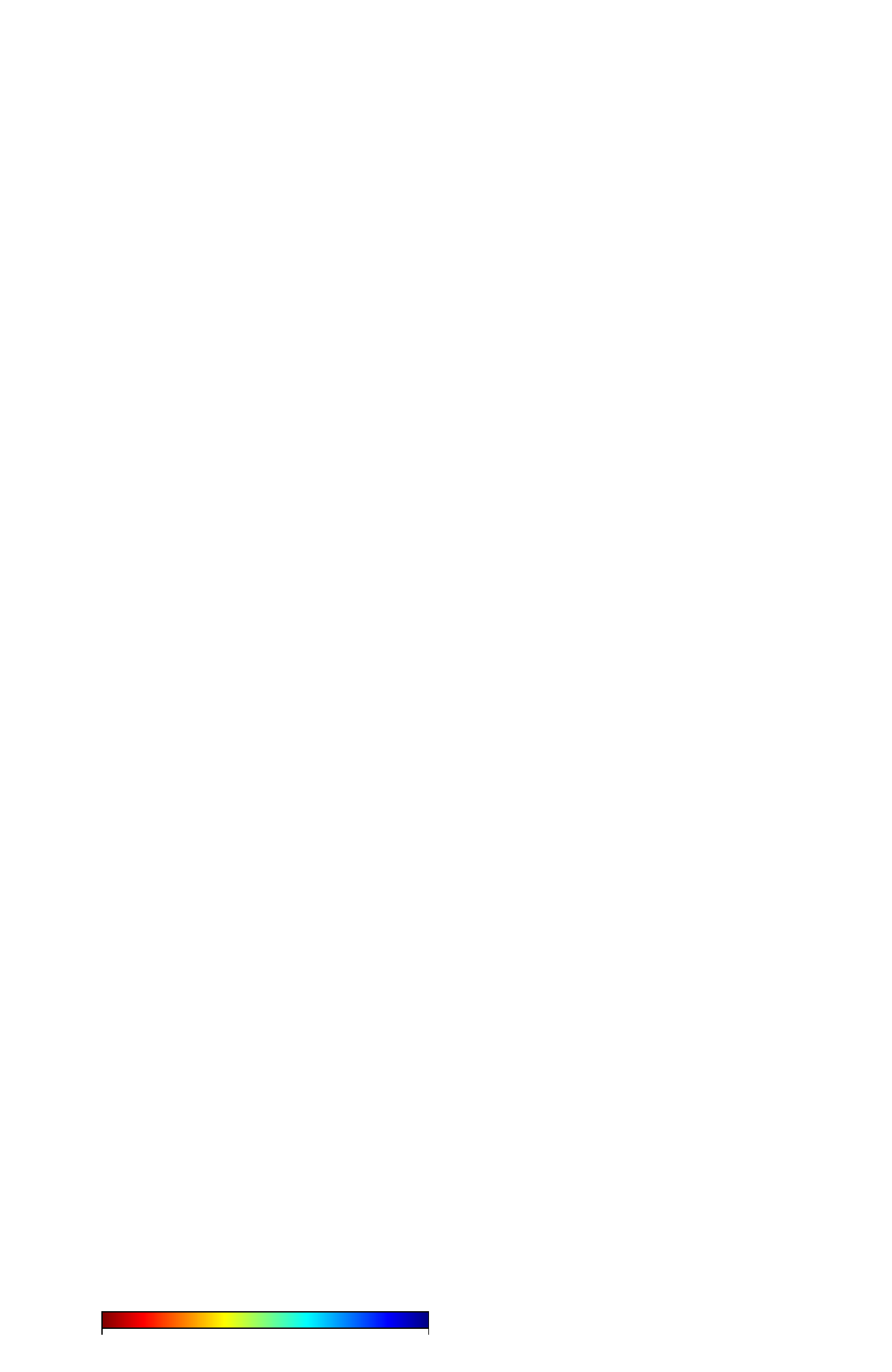
\caption[Estimation of curvature for irregularly sampled torus using different methods.]{Irregular sampling.  } \label{FIG-curv-torus-irreg}
\end{subfigure}
\caption[Estimation of curvature for sampled torus using different methods.]{ The estimation of mean (left column) and Gaussian (right column) curvatures on the sampled torus with 10,000 vertices with different sampling.  }
\end{figure}

As seen in these charts, the best average deviation values of both mean and Gaussian curvature were obtained by DDGO for the sampling of  400 vertices. The average deviations $ H_{avg} $ and $ K_{avg} $ produced by the proposed method AQFC for 3,600 and 10,000 were comparable to the ones given by DDGO. Nevertheless, we consider these values to be satisfactory since the approximation of lower bounds $ H_{min} $, $ K_{min} $ and upper bounds $ H_{max} $, $ K_{max} $ were similar to DDGO and both of them approximate the true bounds of that torus \eqref{EQ-torus-bounds} fair well. The method APSS was not able to approximate upper bound of mean curvature $ H_{max} $ and the lower bound of Gaussian curvature $ K_{min} $ sufficiently. For illustration, we depict the distribution of curvatures on the torus with 10,000 vertices in Figure~\ref{FIG-curv-torus-reg}. The lower and upper bounds are cropped to the values of true bounds. Visually, all three methods depict curvatures correctly, except for the mean curvature from APSS (row (c) left), where the graduation from the areas with low mean curvature to the areas with higher is sharper.

If the torus is finely sampled, that is $ |\mathcal{V}| = \text{360,000}, \text{1,000,000} $, the method AQFC outperformed both DDGO and APSS. The average deviations are significantly smaller for both mean and Gaussian curvature. Also the upper and lower bounds were best preserved by our proposed method.

\subsubsection{Influence of sampling regularity on estimation}
For the irregular sampling of torus and sphere, we calculate the same measures as in the case of regular sampling.  The measured values for both meshes are present in Table~\ref{TABLE-compare-torus-irreg}. 

As indicated by average deviations in this table, the method AQFC is remarkably better in handling of irregular sampling. The approximation of the bounds is worse than in the case of regular sampling of the torus. However, the  irregularities cause major issues for estimation of the bounds for both DDGO and APSS. The distributions of curvatures across the torus are shown in Figure \ref{FIG-curv-torus-irreg}. Note, that the cropping of the color scale to the bounds of true values affects the colouring of the meshes. Still, we are able to notice for how many vertices the curvatures were estimated fair well. As seen in the row (b) of Figure~\ref{FIG-curv-torus-irreg}, the method AQFC produces only minor disruptions in the vertex curvatures compared to the estimation on the regular sampling (see Figure~\ref{FIG-curv-torus-reg} (b)). 

Now, we discuss the case of the sphere, whose irregularities are caused only by flipping some of the diagonals of the regularly sampled quadrilaterals. The depiction of curvatures for the sphere are provided in Figure~\ref{FIG-curv-sphere-irreg}. The color scale is adjusted, so it covers the majority of deviations produced by the tested methods. As seen in the row (d), the problems with estimation by APSS were mainly at the poles of the spheres, which were vertices with high valence. The method DDGO (row (c)) produced bad estimates in the areas where the regularity of the sampling was violated by flipped edges and these deviations are magnified nearer to the poles. The estimation by our method (row (b)) is not visibly affected by the slight irregularities in sampling.

\begin{table}[ht]

\centering
\def\arraystretch{1.6}
\setlength\tabcolsep{1.4em}
\begin{tabular}{|c|c|c|c|}
\hline
& \multicolumn{3}{c|}{Torus} \\
\cline{2-4}
& \makecell{\textbf{AQFC} \\ \textbf{(proposed)}} & DDGO & APSS\\
\hline 
$ H_{min} $ &  -0.627818 & -30.7685 &  -0.57315 \\
\hline
$ H_{max} $ & -0.244277 & 21.7019 &  0.0500141  \\
\hline
$ H_{avg} $ & \textbf{0.00651111} &  0.0874593 & 0.220647 \\
\hline
\hline 
$ K_{min} $ & -0.503117 &  -318.85 &  -1.29606 \\
\hline
$ K_{max} $ & 0.251109 & 3232.91 &   0.925429 \\
\hline
$ K_{avg} $ & \textbf{0.00687775} & 0.909176 & 0.0962692 \\
\hline
\hline
& \multicolumn{3}{c|}{Sphere} \\
\cline{2-4}
&  \makecell{\textbf{AQFC} \\ \textbf{(proposed)}} & DDGO & APSS \\
\hline 
$ H_{min} $  & -1.03639 & -1.14966 & -1.07738   \\
\hline
$ H_{max} $  & -0.997533 & -0.890629 & 0.0058888   \\
\hline
$ H_{avg} $  & \textbf{0.00191217} & 0.0113342 & 0.0214852   \\
\hline
\hline 
$ K_{min} $ & 0.995068 & 0.894789 & -0.00106665    \\
\hline
$ K_{max} $ & 1.0741 & 1.15526 &  1.16044  \\
\hline
$ K_{avg} $ & \textbf{0.00383101} & 0.0163219 & 0.0405045   \\
\hline
\end{tabular}
\caption[Comparison of minimal, maximal and average Gaussian curvature for torus and sphere with irregular sampling.]{The comparison of different methods of the minimal, maximal and average values for estimation of mean and Gaussian curvature on the irregularly sampled torus with 10,000 vertices and sphere with 482 vertices. The lowest average of deviations are highlighted by bold. Note, that the best approximation of minimal and maximal values of curvatures is achieved by using AQFC.}\label{TABLE-compare-torus-irreg}
\end{table}

\begin{figure}[ht]
\small
\centering
\begin{subfigure}{0.48\textwidth}
\def\svgwidth{\textwidth}
\begingroup%
  \makeatletter%
  \providecommand\color[2][]{%
    \errmessage{(Inkscape) Color is used for the text in Inkscape, but the package 'color.sty' is not loaded}%
    \renewcommand\color[2][]{}%
  }%
  \providecommand\transparent[1]{%
    \errmessage{(Inkscape) Transparency is used (non-zero) for the text in Inkscape, but the package 'transparent.sty' is not loaded}%
    \renewcommand\transparent[1]{}%
  }%
  \providecommand\rotatebox[2]{#2}%
  \newcommand*\fsize{\dimexpr\f@size pt\relax}%
  \newcommand*\lineheight[1]{\fontsize{\fsize}{#1\fsize}\selectfont}%
  \ifx\svgwidth\undefined%
    \setlength{\unitlength}{386.94766518bp}%
    \ifx\svgscale\undefined%
      \relax%
    \else%
      \setlength{\unitlength}{\unitlength * \real{\svgscale}}%
    \fi%
  \else%
    \setlength{\unitlength}{\svgwidth}%
  \fi%
  \global\let\svgwidth\undefined%
  \global\let\svgscale\undefined%
  \makeatother%
  \begin{picture}(1,1.65203117)%
    \lineheight{1}%
    \setlength\tabcolsep{0pt}%
    \put(0,0){\includegraphics[width=\unitlength,page=1]{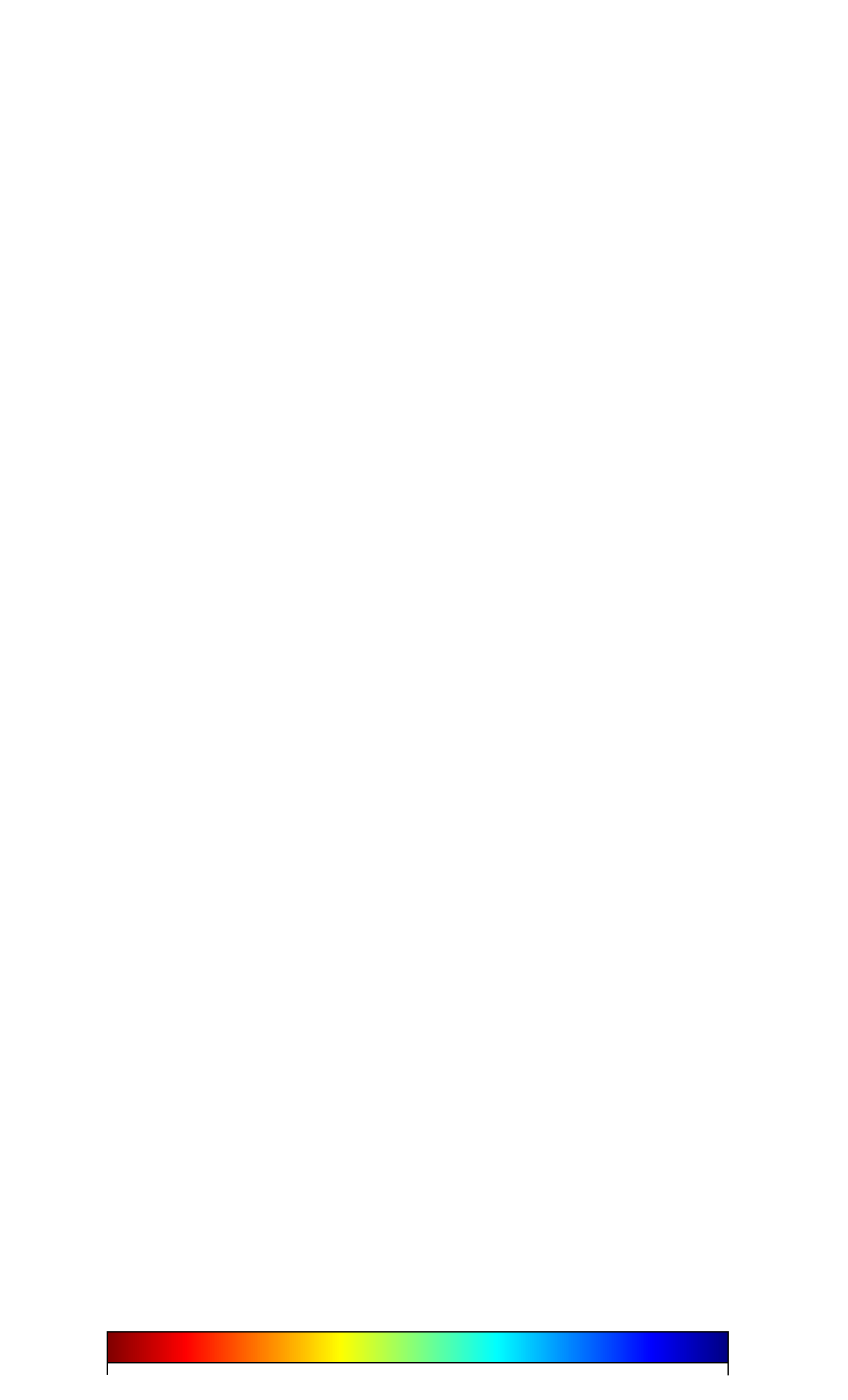}}%
    \put(0.12647503,0.00424246){\color[rgb]{0,0,0}\makebox(0,0)[t]{\lineheight{1.25}\smash{\begin{tabular}[t]{c}$ R = 0 $\end{tabular}}}}%
    \put(0.85904068,0.00424246){\color[rgb]{0,0,0}\makebox(0,0)[t]{\lineheight{1.25}\smash{\begin{tabular}[t]{c}$ R = 200 $\end{tabular}}}}%
    \put(0.22627859,0.89003961){\color[rgb]{0,0,0}\makebox(0,0)[t]{\lineheight{1.25}\smash{\begin{tabular}[t]{c}(a) Input\end{tabular}}}}%
    \put(0.77290465,0.89003961){\color[rgb]{0,0,0}\makebox(0,0)[t]{\lineheight{1.25}\smash{\begin{tabular}[t]{c}(b) AQFC (proposed)\end{tabular}}}}%
    \put(0.22627859,0.11273963){\color[rgb]{0,0,0}\makebox(0,0)[t]{\lineheight{1.25}\smash{\begin{tabular}[t]{c}(c) DDGO\end{tabular}}}}%
    \put(0.77290465,0.11273963){\color[rgb]{0,0,0}\makebox(0,0)[t]{\lineheight{1.25}\smash{\begin{tabular}[t]{c}(d) APSS\end{tabular}}}}%
    \put(0,0){\includegraphics[width=\unitlength,page=2]{curvedness-retheur.pdf}}%
  \end{picture}%
\endgroup%

\caption[Visualization of curvedness on model of Retheur.]{ Curvedness $ R $ on the model of Retheur with  250,562  vertices (a). Curvature estimation by the proposed method AQFC (b) generates less noise than the methods DDGO (c) and APSS (d). } \label{FIG-curvedness-retheur}
\end{subfigure}\hfill
\begin{subfigure}{0.48\textwidth}
\def\svgwidth{\textwidth}
\begingroup%
  \makeatletter%
  \providecommand\color[2][]{%
    \errmessage{(Inkscape) Color is used for the text in Inkscape, but the package 'color.sty' is not loaded}%
    \renewcommand\color[2][]{}%
  }%
  \providecommand\transparent[1]{%
    \errmessage{(Inkscape) Transparency is used (non-zero) for the text in Inkscape, but the package 'transparent.sty' is not loaded}%
    \renewcommand\transparent[1]{}%
  }%
  \providecommand\rotatebox[2]{#2}%
  \newcommand*\fsize{\dimexpr\f@size pt\relax}%
  \newcommand*\lineheight[1]{\fontsize{\fsize}{#1\fsize}\selectfont}%
  \ifx\svgwidth\undefined%
    \setlength{\unitlength}{393.02623508bp}%
    \ifx\svgscale\undefined%
      \relax%
    \else%
      \setlength{\unitlength}{\unitlength * \real{\svgscale}}%
    \fi%
  \else%
    \setlength{\unitlength}{\svgwidth}%
  \fi%
  \global\let\svgwidth\undefined%
  \global\let\svgscale\undefined%
  \makeatother%
  \begin{picture}(1,1.59999568)%
    \lineheight{1}%
    \setlength\tabcolsep{0pt}%
    \put(0,0){\includegraphics[width=\unitlength,page=1]{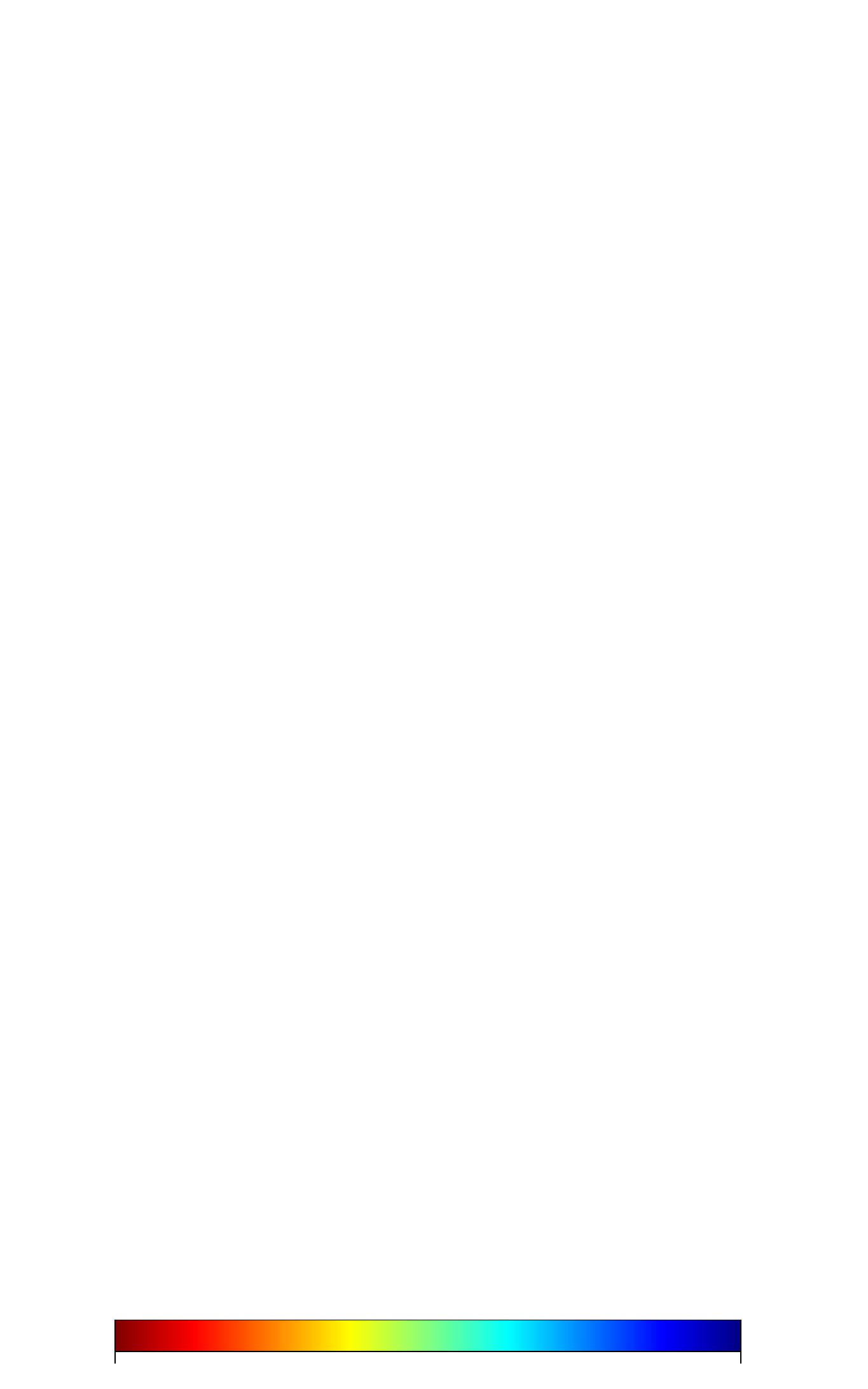}}%
    \put(0.13245445,0.00417684){\color[rgb]{0,0,0}\makebox(0,0)[t]{\lineheight{1.25}\smash{\begin{tabular}[t]{c}$ R = 0 $\end{tabular}}}}%
    \put(0.85369019,0.00417684){\color[rgb]{0,0,0}\makebox(0,0)[t]{\lineheight{1.25}\smash{\begin{tabular}[t]{c}$ R = 200 $\end{tabular}}}}%
    \put(0.23071445,0.8762742){\color[rgb]{0,0,0}\makebox(0,0)[t]{\lineheight{1.25}\smash{\begin{tabular}[t]{c}(a) Input\end{tabular}}}}%
    \put(0.76888635,0.8762742){\color[rgb]{0,0,0}\makebox(0,0)[t]{\lineheight{1.25}\smash{\begin{tabular}[t]{c}(b) AQFC (proposed)\end{tabular}}}}%
    \put(0.23071445,0.11099599){\color[rgb]{0,0,0}\makebox(0,0)[t]{\lineheight{1.25}\smash{\begin{tabular}[t]{c}(c) DDGO\end{tabular}}}}%
    \put(0.76888635,0.11099599){\color[rgb]{0,0,0}\makebox(0,0)[t]{\lineheight{1.25}\smash{\begin{tabular}[t]{c}(d) APSS\end{tabular}}}}%
    \put(0,0){\includegraphics[width=\unitlength,page=2]{curvedness-dame.pdf}}%
  \end{picture}%
\endgroup%

\caption[Visualization of curvedness on model of Dame Assise.]{ Curvedness $ R $ on the model of Dame Assise with  257,818  vertices (a). Curvature estimation by the proposed method AQFC (b) highlights the features (wrinkles, sharp edges) without significant disruption compared to DDGO (c) and APSS (d). } \label{FIG-curvedness-dame}
\end{subfigure}
\caption{Curvedness distribution over meshes.}
\end{figure}

\subsubsection{Curvedness distribution over meshes}
Finally, we compared the tested method on real data meshes obtained from scanning. The comparison of curvedness distribution for Retheus is depicted in Figure~\ref{FIG-curvedness-retheur} and for Dame Assise in Figure~\ref{FIG-curvedness-dame}. Since the curvedness generally has no upper bound, we set the colour scale to the interval $ [0, 200] $, so all features of the mesh are visible.

In the case of Retheus, the noise in the curvedness distribution is present for both DDGO and APSS (Figure~\ref{FIG-curvedness-retheur} (c) and (d)). The most notable disruption is in the areas of the forehead and the left part of the chest, which may be perceived as homogeneous areas in terms of the curvedness. On contrary, the proposed AQFC does not produce any significant noise and the main features, characterized by high curvedness (blue color in the figure), are easily recognizable, e.\ g.\ curls on the beard, creases on the toga or the sharp edges of the pedestal. The slight deviations of curvatures are present in the smooth areas, e.\ g.\ chest, however these are caused by the material of the scanned statue and are considered the feature of the mesh as well.

The effect of imprecise curvedness estimation by DDGO and APSS is even more visible in the case of Dame Assise, see Figure~\ref{FIG-curvedness-dame} (c) and (d). The features like creases on the dress or parts of the face are hardy distinguishable. However, AQFC captures even the fine details as waviness of the bottom part of the dress (Figure~\ref{FIG-curvedness-dame} (b)).

\section{Conclusion}\label{SEC-conclusion}
In our work, we introduced mean and Gaussian curvature estimation method AQFC, which is based on local algebraic fitting of quadratic surfaces. AQFC is capable to process arbitrary polygon meshes, whose vertices are equipped with associated normal vectors.

 The experimental results indicate, that AQFC is feasible for estimation of the mean and Gaussian curvature and their derivatives. We consider our method to be more universal albeit the well established method DDGO provides more precise results for the regular sampling of the surfaces. However, DDGO fails to provide satisfactory results for meshes with irregular sampling. Even if APSS is less prone to irregularities compared to DDFO, AQFC provides more precise estimate of curvatures. 
 
The proposed method AQFC is also robust at the computation of curvedness for meshes, where the prescription of the underlying surface is unknown. Thus, it has potential to be used in various application as denoising or smoothing of the meshes.

In future work we plan on to extend this approach for processing of point clouds. The modifications that need to be made is the way of obtaining neighbourhood since edges and faces are not present. This may be solved by using the feasible space partitioning structure as octrees or $ k $-d trees.

\bibliographystyle{eg-alpha-doi}

\bibliography{qfc-bib}

\end{document}